\documentclass[3p,twocolumn,preprint]{elsarticle}
\usepackage{amsmath}
\usepackage{verbatim}
\usepackage{amssymb,bbm}
\usepackage{latexsym,url}

\newcommand{\oR}{{\mathbb R}}
\newcommand{\oN}{{\mathbb N}}
\newcommand{\oZ}{{\mathbb Z}}
\newcommand{\oQ}{{\mathbb Q}}

\newcommand{\EE}{\mathcal E}
\newcommand{\RR}{{\mathcal R}}

\newcommand{\CUT}{\text{\rm CUT}}

\newcommand{\METO}{\text{\rm MET}^{01}}
\newcommand{\METU}{\text{\rm MET}}
\newcommand{\minor}{\preccurlyeq}
\newcommand{\psd}{\succeq}
\newcommand{\cosp}{\cos\frac{\pi}{p}}
\newcommand{\SSS}{{\mathcal S}}

\newcommand{\GG}{{\mathcal G}}

\newcommand{\mc}{\text{\rm mc}}
\newcommand{\ignore}[1]{}

\newcommand{\arcsinh}{\text{\rm arcsinh}}

\newcommand{\ka}{\kappa}

\newcommand{\sdp}{\text{\rm sdp}}
\newcommand{\ip}{\text{\rm ip}}
\newcommand{\rank}{\text{\rm rank}}

\newtheorem{theorem}{Theorem}[section]

\newtheorem{lemma}{Lemma}[section]

\newproof{pf}{Proof}
\begin{document}
\title{Computing the Grothendieck constant of some graph classes}

\author[cwi,til]{M.~Laurent}
\ead{M.Laurent@cwi.nl}

\author[cwi]{A.~Varvitsiotis\corref{cor1}}
\ead{A.Varvitsiotis@cwi.nl}

\cortext[cor1]{Corresponding author: CWI, Postbus 94079,
	      1090 GB Amsterdam. Tel: +31 20 5924170; Fax: +31 20 5924199.}

\address[cwi]{Centrum Wiskunde \& Informatica (CWI), Science Park 123,
	        1098 XG Amsterdam,
	        The Netherlands.}
\address[til]{Tilburg University, 
P.O. Box 90153, 
5000 LE Tilburg, 
The Netherlands.} 

\begin{abstract}
Given a graph $G=([n],E)$ and  $w\in\oR^E$, consider the integer program  
${\max}_{x\in \{\pm 1\}^n} \sum_{ij \in E} w_{ij}x_ix_j$ 
and its canonical semidefinite programming relaxation 
${\max} \sum_{ij \in E} w_{ij}v_i^Tv_j$, where the maximum is taken 
over all unit vectors $v_i\in\oR^n$. The integrality gap of this relaxation is known as the Grothendieck constant $\ka(G)$ of $G$.
We present  a closed-form formula for the Grothendieck constant of $K_5$-minor free graphs  and derive that it is at most $3/2$. Moreover, we   show that $\ka(G)\le \ka(K_k)$ if the cut polytope of $G$ 
 is defined by inequalities supported by at most $k$ points.
Lastly, since the Grothendieck constant of $K_n$ grows as $\Theta(\log n)$,  it is interesting  to identify  instances with large gap. 
However this is not the case for the clique-web inequalities, a wide class of valid inequalities for the cut polytope, whose integrality ratio is shown to be bounded by 3.

\end{abstract}

\begin{keyword} Grothendieck constant \sep elliptope \sep  cut polytope \sep clique-web inequality
\end{keyword}
\maketitle

\section{Introduction}

Let $G=([n],E)$ be a (simple loopless)  graph and  $w=(w_{ij}) \in \oR^E$. Consider the integer quadratic program  over the hypercube

\begin{equation}\label{defip} 
\ip(G,w):= \max_{x\in\{\pm 1\}^n}\sum_{ij\in E}w_{ij}   x_ix_j,
\end{equation}
and its   canonical  semidefinite programming relaxation
\begin{equation}\label{eqG1}
\sdp(G,w):=\displaystyle \max_ {u_1,\ldots,u_n\in \oR^n,\ \|u_i\|=1} 
 \sum_ {ij \in E}w_{ij}u_i^Tu_j.
\end{equation}
Let $\ka(G)$ denote the integrality gap of relaxation (\ref{eqG1}), defined by
\begin{equation}\label{ig}\ka(G)=\sup_{w \in \oR^E} \frac{\sdp(G,w)}{\ip(G,w)}.
\end{equation}
In other words, $\ka(G)$ is  the smallest constant $K>0$ for which 
$\sdp(G,w)\le K \cdot  \ip(G,w), \ \forall w \in \oR^E.
$  Alon et al. \cite{AMMN} call this graph parameter  the  {\em Grothendieck constant}  of $G$ and
 prove that  
\begin{equation}\label{eqAMMN}
\Omega(\log \omega(G))=\ka(G)=O(\log \vartheta (\bar G)).
\end{equation}
Here $\omega(G)$ denotes the  maximum size of a clique in $G$ and $\vartheta(\bar G)$ the Lov\'asz theta function of the complementary graph $\bar{G}$, for which it is known  that   $\omega(G) \le \vartheta(\bar G) \le  \chi(G)$   \cite{lov79}.  
Hence, for the complete graph $G=K_n$, 
$\ka(K_n)=\Theta(\log n)$.

The name of the constant goes back to   Grothendieck \cite{Gro53}, who considered the case of bipartite graphs and showed the  existence  of  a constant $K>0$ for which    
$\label{eqG}
 \sdp(K_{m,n}, w)  \le K \cdot \ip(K_{m,n}, w)$ for all $ m,n\in \oN$ and $ w\in \oR^{mn}.
$ The smallest such constant is known as  {\em Grothendieck's constant} and is denoted by $K_G$.  It is a long standing open problem to compute the exact value of $K_G$. It is known that $K_G < {\pi [2\ln (1+\sqrt 2)]^{-1}} \sim 1.782 
$ \cite{Kri77,BMMN}, and that $K_G \ge 1.6769...$ \cite{Reeds}. 
Recently  Bri\"et et al. \cite{BOV} show that $\ka(G)\le {\pi\over 2\ \arcsinh ((\vartheta(\bar{G})-1)^{-1})}$, which gives the above bound $1.782$ for bipartite graphs and improves the upper bound  in (\ref{eqAMMN}) when $\vartheta(\bar{G})$ is small.

In recent years, Grothendieck type inequalities have  received a significant  amount of attention due to their  various applications, most notably in the design of approximation algorithms  and quantum information theory (see, e.g., \cite{AMMN,AN04,BMMN,FR94,Pi08,W06}). 

The paper is organized as follows. In Section \ref{basic} we collect  basic properties of $\ka(G)$. In Section \ref{computing} we establish a closed-form  formula for the Grothendieck constant  of  $K_5$-minor free graphs  in terms of their girth and bound $\ka(G)$ in terms of the size of the supports of the facets of the cut polytope.  In Section \ref{igap} we show that the integrality gap achieved by the  clique-web inequalities, a wide class of valid inequalities for  the cut-polytope, is bounded by 3.

\section{Basic properties}\label{basic}
We first introduce some notation. Throughout  $[n]=\{1,\ldots,n\}$.
 Let  $\SSS_n^+$  denote the cone of positive semidefinite matrices; the notation $A\succeq 0$ means that $A\in \SSS_n^+$.
For matrices $A,B$, $\langle A,B\rangle =\sum_{i,j} A_{ij} B_{ij} $ stands for the usual trace inner product.
Let $e$ denote the all ones  vector and $J=ee^T$ the all ones matrix, of the appropriate dimension. 
 
Let  
 $\EE_n:=\{X\in \SSS_n^+\mid X_{ii}=1\ \forall i\in [n]\}$ and 
 $\CUT_n :={\rm conv}(X\in\EE_n\mid \rank \ X=1\}.$
 Moreover, define
 $ \EE(G):=\pi_E(\EE_n),\ \CUT(G):=\pi_E(\CUT_n),$
where $\pi_E$ denotes the projection from $\oR^{n\times n}$ onto  the subspace  $\oR^E$
 indexed by the edge set of $G$. They  are known, respectively,  as the {\em elliptope} and the {\em cut polytope} of $G$ and satisfy   $\CUT(G)\subseteq \EE(G)$.  We refer, e.g., to \cite{DL97} and further references therein for a detailed study of these geometric objects.
 
 For $w\in \oR^E$,  let   $\ka(G,w)=\sdp(G,w)/\ip(G,w)$.

\subsection{A geometric reformulation for $\ka(G)$} 
Clearly,  the Grothendieck constant $\ka(G)$ is the smallest  dilation of $\CUT(G)$ containing $\EE(G)$.
 
\begin{lemma}\label{lemKGdilation} 
For any graph $G$,  \begin{center}$ \ka(G) = \min \{ K  \ | \  \EE(G)\subseteq K \cdot  \CUT(G)\}.$\end{center}
\end{lemma}
\begin{pf}
Directly, since  $\ip(G,w)=\underset{x \in \CUT(G)}{\max}w^Tx \ $ and $\ \sdp(G,w)=\underset{x \in \EE(G)}{\max}w^Tx.$
\qed\end{pf}

As the origin lies in the interior of  $\CUT(G)$,  the polytope $\CUT(G)$ has a linear inequality description consisting of  finitely many facet-defining inequalities of the form $w^Tx\le 1$. Let us recall the following {\em switching operation:}
Given $w\in\oR^E$, its switching by  $S\subseteq [n]$ is the vector $w^{(S)} \in \oR^E$ whose $(i,j)$-th entry  is $-w_{ij}$ if the edge $ij$ is cut by the partition $(S,[n]\setminus S)$ and $w_{ij}$ otherwise.
It is well known that the switching operation preserves valid inequalities and facet defining inequalities of the cut polytope 
\cite{BM,DL97}. Moreover, $\sdp(G,w)=\sdp(G,w^{(S)})$ and $\ip(G,w)=\ip(G,w^{(S)})$. 
This implies the next lemma which gives a useful  reformulation for  $\ka(G)$.
\begin{lemma}\label{ratio}For any graph $G$, $$\ka(G)=\sup_{w \in \oR^E} \ka(G,w),$$ where the supremum ranges over  all facet defining inequalities of $\CUT(G)$,  distinct up to switching.
\end{lemma}

   \subsection{Connections  with   max-cut}

The study of the cut polytope $\CUT(G)$ and of the elliptope $\EE(G)$ is largely motivated by their relevance to the  maximum cut problem in combinatorial optimization.
Given $G=([n],E)$ and  $w \in \oR^E$, the max-cut problem asks for a cut 
 of maximum weight. 
Thus we want to compute
$\mc(G,w)= \max_{x\in \{\pm 1\}^n} 
{1\over 2}\sum_{ij\in E} w_{ij}(1-x_ix_j)=\max_{x\in \{\pm 1\}^n}  \frac{1}{4}x^TL_{{G},w}x.$ Here,  $L_{{G},w}$ is the Laplacian matrix,
with $(i,i)$-th entry $\sum_jw_{ij}$ and with $(i,j)$-th entry $-w_{ij}$ if $ij\in E$ and 0 otherwise. 
The canonical semidefinite programming relaxation of max-cut (considered e.g.~in \cite{GW95}) is 
$\sdp_{\text{GW}}(G,w)=\max_{X\in \EE_n} {1\over 4} \langle L_{G,w},X\rangle$. 
Hence the quadratic integer problem (\ref{defip}) and the max-cut problem are affine transforms of each other, and the same for their canonical semidefinite relaxations; namely,
$\mc(G,w)= {1\over 2}\left(w(E)+\ip(G,-w)\right)$
and
$\sdp_{\text{GW}}(G,w)= {1\over 2}\left(w(E)+\sdp(G,-w)\right).$

In particular, this implies that, given $w\in\oQ^E$, deciding whether $\ip(G,w)=\sdp(G,w)$ is an NP-complete problem   \cite{DP93c}. 

The following lemma is easy to verify.
\begin{lemma}\label{lem} Let  $A \in \SSS_n^+$ and 
$B=\left(\begin{matrix} 0 & A/2 \cr A/2 & 0\end{matrix}\right)$.
Then, $\displaystyle \max_{Z\in \EE_{2n}} \langle B,Z\rangle= \max_{X\in \EE_n} \langle A,X\rangle$, \\ and 
$\displaystyle\max_{z\in \{\pm 1\}^{2n}} z^TBz =\max_{x\in \{\pm 1\}^n} x^TAx$.
\end{lemma}

When $w \ge 0$,  $L_{\text{G},w} \psd 0$ and thus Lemma \ref{lem} implies that 
$\sdp_{\text{GW}}(G,w) 
=\underset{Z \in \EE_{2n}}{\max} \langle B,Z\rangle$, where $B$ is as in the lemma with $A/2 := L_{G,w}/8$. By the definition of the Grothendieck constant $K_G$, this implies that
$\sdp_{\text{GW}}(G,w)\le K_G\cdot \mc(G,w)$. However, this approximation guarantee is not interesting since  we know by \cite{GW95} that 
$\sdp_{\text{GW}}(G,w)\le 1.138 \cdot \mc(G,w)$, while $K_G\ge 1.6$.

On the other hand, the Grothendieck constant  $\ka(G)$ bounds the semidefinite approximation for max-cut  for edge weights satisfying $w(E)\ge 0$. 

\begin{lemma}
Let $G=(V,E)$ be a graph and $w\in\oR^E$ with $w(E)\ge 0$ and $\mc(G,w)>0.$  Then, $ \sdp_{\text{GW}}(G,w)\le \ka(G)\cdot \mc(G,w)$.
\end{lemma}
\begin{pf}
Indeed, $\sdp(G,-w)\le \ka(G)\cdot \ip(G,-w)$ and $w(E)\le \ka(G)\cdot w(E)$ imply 
$$
{\sdp_{\text{GW}}(G,w)\over \mc(G,w)} = {w(E)+\sdp(G,-w)\over w(E)+\ip(G,-w)} \le \ka(G).$$
\qed\end{pf}

   \subsection{Behaviour under graph operations}
 It follows immediately from the definition that the graph parameter $\ka(\cdot)$ is monotone nonincreasing with respect to deleting edges. That is,

\begin{lemma}\label{lemsubgraph}
If $H \subseteq G$ then $\ka(H)\le \ka(G)$.
\end{lemma}

This is not true for the operation of contracting an edge. For instance, $\ka(K_2)=1 < \ka(C_3)=3/2$,  while $\ka(C_4)< \ka(C_3)=3/2$ (cf.~Theorem \ref{theocircuit}). 
So $\ka(G)$ and $\ka(G\slash e)$ are not comparable in general.

Given two graphs $G_1=(V_1,E_1)$ and $G_2=(V_2,E_2)$ for which $V_1\cap V_2$ is a clique
in both $G_1$ and $G_2$, the graph $G=(V_1\cup V_2,E_1\cup E_2)$ is called their {\em clique sum}, 
or their  {\em clique $k$-sum} when $|V_1\cap V_2|=k$.

\begin{lemma}\label{lemcliquesum}
Assume $G$ is the clique $k$-sum of $G_1$ and $G_2$,  $k\le 3$.
Then, $\ka(G)=\max(\ka(G_1),\ka(G_2))$.
\end{lemma}

\begin{pf}
Let $\lambda:=\max(\ka(G_1),\ka(G_2))$ and $n=|V|$.
The inequality $\ka(G)\ge \lambda$ follows from  Lemma \ref{lemsubgraph}.

 For the other direction, let   $x\in \EE(G)$ and $X\in \EE_n$ such that $x=\pi_E(X)$; we have to 
 show that $x\in\lambda \cdot \CUT(G)$. Let $X_i$ denote the principal submatrix of $X$ indexed by $V_i$, for $i=1,2$. 
  As $\EE(G_i)\subseteq \ka(G_i)\cdot \CUT(G_i)\subseteq \lambda \cdot\CUT(G_i)$, we deduce that 
  $\pi_{E_i}(X_i) \in \lambda\cdot \CUT(G_i)$. Since the
  linear inequality description of $\CUT(G)$ is obtained by juxtaposing the linear inequality descriptions of $\CUT(G_1)$ and $\CUT(G_2)$ and identifying the variables corresponding to edges contained in $V_1\cap V_2$  \cite{B83}, the claim follows. 
\qed\end{pf}

\section{Computing the Grothendieck constant for some graph classes}\label{computing}
We start this section by introducing the main  objects and some fundamental results associated with them, that form the  basic  ingredients of our approach.  

A graph $H$ is called a {\em minor} of a graph $G$, denoted by $H \minor G$, if $H$ can be obtained from $G$, through a series of  edge deletions and edge contractions. 

Given a graph $G=([n],E)$, consider the {\em metric polytope}  
$\METU(G)\subseteq \oR^E$ defined by the inequalities $-1\le x_e\le 1$ for $e\in E$, and 
$$ 
x(C\setminus F)-x(F)\le |C| -2, $$ 
for every circuit   $C$  in $G$ and  $F\subseteq C$ with $|F|$ odd \cite{DL97}. 
Additionally, define  $\METO(G):=f(\METU(G))$, where $f(x)=e-2x$ for $x\in\oR^E$.  The cut and metric polytopes are related as follows.
\begin{theorem}\cite{BM}\label{theocut}
For any graph $G$, 
$\CUT(G)\subseteq \METU(G)$, with equality if and only if $K_5 \not \minor G$.
\end{theorem}

Any matrix $X\in \EE_n$ has its diagonal entries all equal to 1. Hence all  its  entries lie in  $[-1,1]$ and  can thus be parametrized as $x_{ij}=\cos(\pi y_{ij})$ where $y_{ij}\in [0,1]$. Let $\cos (\pi \METO(G))=$
\begin{center} $\{ (\cos(\pi y_e))_{e\in E} \mid y\in \METO(G)\}.$ \end{center}

\begin{theorem}\cite{Lau97}\label{theocos}
 $\EE(G)\subseteq \cos (\pi \METO(G))$, with equality if and only if $K_4 \not \minor G$.
\end{theorem}

Thus,  equality holds when  $G=C_n$. Moreover,
\begin{lemma}\label{econdition} \cite{BJT93} For $p$ even,  we have 
 $ce^T \in \EE(C_p)$ for all $c \in [-1,1]$.
For  $p$ odd, we have   $ce^T \in \EE(C_p)$ if and only if $-\cos\frac{\pi}{p} \le c \le 1$. 
\end{lemma}

\subsection{The case of circuits}
Using  the parametrizations of $\METU(C_n)$ and $\EE(C_n)$ given by  Theorems \ref{theocut} and \ref{theocos}, respectively, we are able  to compute $\ka(C_n)$. Specifically, 

\begin{theorem}\label{theocircuit}
The Grothendieck constant of a circuit $C_n$ of length $n\ge 3$ is equal to
$$\ka(C_n)= {n\over n-2} \cos\Big({\pi\over n}\Big).$$

\end{theorem}

\begin{pf}
 By Lemma \ref{ratio} it suffices to compute $\ka(C_n,w)$ for facet defining  inequalities of $\CUT(C_n)$. 
 By Theorem \ref{theocut}, they correspond to the circuit  inequalities and, since they are all switching equivalent, it suffices  to
 consider one of them; for instance, we can choose $w^Tx=-x(E) $  for odd $n$, and  $w^Tx=x_e-x(E\setminus \{e\})$   for even $n$. In both cases, we find that $\ip(C_n,w)=n-2$. Thus it now suffices to show that $\sdp(C_n,w)=n\cos(\pi/n)$ as this will give the desired value for $\ka(C_n,w)$.
 
 For $n$ odd, it
 is known that $\sdp_{\text{GW}}(C_n, e)=\frac{n}{4}\left(2+2\cos\frac{\pi}{n}\right)$ (see \cite{DP93b}), which implies that  
 $\sdp(C_n,-e)=2\ \sdp_{\text{GW}}(C_n,e)-n = n\cos(\pi/n)$.  This can also  be easily verified using the parametrization of $\EE(C_n)$ from Theorem \ref{theocos}.

One can also compute $\sdp(C_n,w)$ for $n$ even and $w=(-1,1,\ldots,1)$  using Theorem \ref{theocos}; 
it turns out that this  has also been
computed in \cite{W06} in the context of quantum information theory.
\qed\end{pf}

\subsection{The case of  $K_5$-minor free graphs}\label{secnoK5}

Since $K_5$-minor free graphs are 4-colorable \cite{W37}, we deduce from   (\ref{eqAMMN}) that their Grothendieck constant  $\ka(G)$ is bounded.
Here we give a closed-form formula for $\ka(G)$ in terms of the girth of $G$.

\begin{theorem}\label{theonoK5}
If $G$ is a graph with no $K_5$ minor (and $G$ is not a forest), then

$$\ka(G)= {g\over g-2}\cos\left({\pi\over g}\right),$$
where $g$ is the minimum length of a circuit in $G$. 
\end{theorem}

\begin{pf}
Directly from Theorem \ref{theocircuit} 
using the facts that all facets of $G$ are supported by circuits (Theorem \ref{theocut}) and that the function ${n\over n-2}\cos({\pi\over n})$ is monotone nonincreasing in $n$.
\qed
\end{pf}
As  a direct application we recover the  values   $\ka(K_{2,n})=\ka(K_{3,n})=\sqrt{2},$ for $ n \ge 3$ \cite{FR94}.

\subsection{Graphs whose cut polytope is defined by inequalities supported by at most $k$ points}

We show here that the Grothendieck constant can be bounded in terms of the size of the supports of the inequalities defining facets of the cut polytope.
The {\em support graph} of an inequality $w^Tx \le 1$ is  the graph $H=(W,F)$, where $F=\{ij\in E\mid w_{ij}\ne 0\}$ and $W$ is the set of nodes covered by $F$. We say that  $w^Tx\le 1$ is {\em supported by at most $k$ points}
when 
  $|W|\le k$.
For instance, 
a triangle inequality is supported by three points.

Fix an integer $k\ge 2$. 
Let $\RR_k(K_n) \subseteq \oR^{E_n}$ be  the polyhedron defined by all
valid inequalities for $\CUT_n$ supported by at most $k$ points. For $G=([n],E)$, let   
$\RR_k(G):= \pi_E(\RR_k(K_n))$.
For instance, 
$\RR_3(K_n)=\METU(K_n)$, and thus  $\RR_3(G)=\METU(G)$.

Clearly, $\CUT(G)\subseteq \RR_k(G)$. Define the class $\GG_k$ of all graphs $G$ for which 
$\CUT(G)=\RR_k(G)$. 
For instance, $\GG_2$ consists of all forests (thus the $K_3$-minor free graphs) and $\GG_3$ of the $K_5$-minor free graphs. Thus both are closed under taking minors; this holds for any   $\GG_k$.

\begin{theorem} The class $\GG_k$  is closed under taking minors.
 \end{theorem}
 
 \begin{pf}
 It follows directly from the definition that $\GG_k$ is closed under edge deletion. It remains to verify that it is closed under edge contraction. Let  $G=(V,E)$ and $G':=G\slash e=(V',E')$, where $e=(1,2)$ and $V^{'}=\{2,\ldots,n\}$.

 Given $y\in\oR^{E'}$, define its  extension $\tilde{y}\in \oR^E$ by
 $\tilde{y}_{12}=1$, $\tilde{y}_{1i}=y_{2i}$ if $1i\in E$ with $i\ge 3$, $\tilde{y}_{2i}=y_{2i}$ if $2i\in E$ with $i\ge 3$, and $\tilde{y}_{ij}=y_{ij}$ if $ij\in E$ with $i,j\ge 3$.
 One can easily verify that $\tilde{y}\in\CUT(G)$ iff $y\in\CUT(G\slash e)$.
 
 We now  verify that if $y\in \RR_k(K_{n-1})$, then  $\tilde{y} \in \RR_k(K_n)$.
 Indeed, say $w^Tx\le 1$ is a valid inequality for $\CUT_n$ supported by  at most $k$ points. 
 Define the inequality on $x=(x_{ij})_{2\le i<j\le n}$:
 $$b^Tx:=\sum_{i=3}^n (w_{1i}+w_{2i})x_{2i} + \sum_{3\le i<j\le n} w_{ij} x_{ij} \le 1-w_{12}.$$
 Obviously it is supported by  at most $k$ points and it is valid for $\CUT_{n-1}$.
 Hence $b^Ty\le 1-w_{12}$, which implies $w^T\tilde{y} \le 1$. 
 
 Assume $G\in \GG_k$, i.e.,
 $\CUT(G)=\pi_E(\RR_k(K_n))$ and let $z\in \pi_{E'}(\RR_k(K_{n-1}))$; we show that $z\in \CUT(G\slash e)$. Say $z=\pi_{E'}(y)$ where $y\in \RR_k(K_{n-1})$. By the discussion above, the extension $\tilde{y}$ of $y$ belongs to $\RR_k(K_n)$ and thus $\pi_E(y)\in \pi_E(\RR_k(K_n))=\CUT(G)$.
This in turn implies that $z\in \CUT(G\slash e)$ since $\pi_E(y)$ is the extension of $z$.
\qed\end{pf}


Clearly,  for  $G \in \GG_2$, $\ka(G) = \ka(K_2)=1$. Moreover, Theorem \ref{theonoK5} implies  that  $\ka(G)\le \ka(K_3)=3/2$ for $G\in\GG_3$.  This pattern extends to any $k$.

\begin{theorem}
If $G\in \GG_k$ then $\ka(G)\le \ka(K_k)$.  Moreover, this bound is tight since $K_k \in \GG_k$.
\end{theorem}

\begin{pf}
It is enough to show that for any graph $G$, $\EE(G)\subseteq \ka(K_k) \cdot \RR_k(G)$. Moreover,  it suffices to consider only  $G=K_n$, as the general result follows by taking projections. 

Let $y\in \EE(K_n)$ and let $w^Tx\le 1$ be a valid inequality for $\CUT_n$  with support $H=(W,F)$ where 
$|W|\le  k$. 
Then, 
$w^Ty =\pi_F(w)^T\pi_F(y) \le \sdp(H,\pi_F(w)) \le \ka(H)\cdot \ip(H,\pi_F(w)) \le \ka(K_k),$
where we use the facts that $\ka(H)\le \ka(K_k)$ and $\ip(H,\pi_F(w))\le 1$ for the right most inequality.
\qed\end{pf}

One can verify that $\ka(K_7)=3/2$ (see \cite{Lau04}). Hence, $\ka(G)\le 3/2$ for all $G\in \GG_7$.


 \section{Integrality gap of clique-web inequalities}\label{igap}
 We have already seen that $\ka(K_n)=\Theta(\log n)$ and it is an interesting question to identify explicit instances that achieve this integrality gap. This was posed as an open question in \cite{AMMN} and  instances with large gap are given in \cite{AB05}.
 In this section we show that the gap is bounded by 3 for clique-web inequalities, a wide class of valid inequalities for $\CUT_n$.

Given integers $p$ and $r$ with $p \ge 2r+3$, the {\em antiweb graph} $\text{AW}^r_p$ is the graph with vertex set $[p]$,  and with edges $(i,i+1),\ldots,(i,i+r)$  for $i \in [p]$, where the  indices are taken modulo $p$. The {\em web graph}  $\text{W}^r_p$ is defined as  the complement of $\text{AW}^r_p$ in $K_p$. 
Call the set of edges $(i,i+s)$ for $i\in [p]$ (indices taken modulo $p$) the {\em $s$-th band}, so that AW$^r_p$ consists of the first $r$ bands and 
W$^r_p$ consists of the last $\lceil q/2\rceil$ bands in $K_p$.

Let $p, q,r,n$ be integers satisfying $p-q=2r+1,$ $q\ge 2, n=p+q$. The (pure) {\em clique-web inequality} with parameters $n, p,q,r$ is the inequality 
\begin{equation}\label{cweb}
-x(K_q)-\sum_{1\le i \le q \atop{q+1 \le j \le n}}x_{ij}-x(\text{W}^r_p)\le q(r+1).
 \end{equation}
 The support graph of (\ref{cweb}), denoted by CW$^r_p$,  consists of  a clique on the first $q$ nodes,  a web on the last $p$ ones, and a complete bipartite graph between them. 
 It is known that clique-web inequalities  define facets of $\CUT_n$.  Note that hypermetric and bicycle odd wheel inequalities arise as special cases of (\ref{cweb}), for $r=0$  and $r=\frac{n-5}{2}$, respectively (see \cite{DL97}). 
 
Since the left-hand side of (\ref{cweb}) is equal to    $-x(K_{p+q}) +x( \text{AW}^r_p)$,  one easily obtains that 
$\sdp(\text{\rm CW}_p^r,-e) \le (p+q)/2 +pr =  q(r+1) + (2r+1)^2/2$, which implies that 
 $\ka(\text{\rm CW}_p^r,-e)\le 1+\frac{(2r+1)^2}{q(2r+2)}.$
This directly implies the following:

\begin{lemma}\label{trbound}
The integrality gap  of a  clique-web inequality with $q \ge 2r+1$  is upper bounded by 2. 
\end{lemma}

We now consider the case when $q\le 2r$.
\begin{theorem}\label{cwthm} 
The integrality gap of a clique-web inequality with $q\le 2r$ is upper bounded by 3.
\end{theorem}

\begin{pf}  
We can rewrite $\sdp(\text{\rm CW}_p^r,-e)$ as
\begin{equation}\label{sdp1} 
 \underset{X \in \EE_n}{\max}   -\sum_{ij \in K_q}X_{ij} - \sum_{1\le i \le q \atop{q+1 \le j \le n}}X_{ij} -\sum_{ij \in \text{W}^r_p} X_{ij}.
\end{equation}
Notice that the program (\ref{sdp1}) is invariant under the action 
of the full symmetric group $S_q$ acting on the row/column indices in $[q]$. 
Moreover,  (\ref{sdp1})  is invariant under the action  of the group of cyclic permutations in $S_p$ acting on the row/column indices in $\{q+1,..,n\}$.  Thus, we can restrict without loss of generality the matrix $X$ in (\ref{sdp1}) to satisfy the following invariance conditions:
\begin{equation*}\label{inv}
\begin{array}{ll}
X_{ij}=a & \text{ for } 1\le i\ne j\le q,\\
X_{ij}=b & \text{ for } 1\le i\le q< j\le n,\\
X_{ij} = c_{|j-i \hspace*{-1mm} \mod p|} & \text{ for } q+1\le i\ne j\le n
\end{array}\end{equation*}
for some scalars $a,b,c_1,\ldots,c_{\lceil q/2\rceil}$.
Hence $X$ has the form
$\small{X:=\left(\begin{array}{cc}
aJ_{q,q}+ (1-a)I_{q} & bJ_{q,p} \\
bJ_{p,q} & X_p
\end{array}\right)},$
where $X_p$ denotes the principal submatrix of $X$ indexed by $\{q+1,\ldots,n\}$. 
One can easily verify that $X\succeq 0$ if and only if 
 $ Y:=\small{\left(\begin{array}{cc}
\beta & be^T \\
be & X_p
\end{array}\right)} \succeq 0,$ after setting $\beta:= \frac{(q-1)a+1}{q}$.

Consider first the case when $q$ is even; so $p$ is odd, all bands in W$^r_p$ have size $p$, and   the objective function in (\ref{sdp1}) reads
\begin{equation}\label{obj1}
 {q\over 2}(1-q\beta) -pq b - p(c_1+\ldots + c_{q/2}).
 \end{equation}
 If $\beta=0$, then $b=0$ and  Lemma \ref{econdition} implies  that $c_s\ge -\cos(\pi/p)$ for all $s$.
Indeed each band of W$^r_p$ is a circuit or a disjoint union of circuits 
(e.g. the first band of $\text{W}^2_9$ is a union of three triangles). As $p$ is odd, at least one of these circuits is an odd circuit of size $p'\le p$, so that  Lemma \ref{econdition} implies that the entries on the band are at least $-\cos(\pi/p')\ge -\cos(\pi/p)$. Now the objective value is equal to ${q\over 2}-p(c_1+\ldots +c_{q/2})\le {q\over 2}(p\gamma+1)\le {q\over 2}(p+1)= {q\over 2}(q+2r+2)\le 2 \ q(r+1)$ (as $q\le 2r$), setting $\gamma:=\cos(\pi/p)$.

Assume now  $\beta>0$. Taking the Schur complement in the above matrix $Y$ with respect to the entry $\beta$, we can rewrite the condition $ Y \succeq 0$ as $X_p-{b^2\over \beta} J\succeq 0$.
If $\beta=b^2$, then $c_s=1$ for all $s$ and the maximum of ${q\over 2}(1-qb^2)-pqb-pq/2$ for $b\in [-1,1]$ is easily verified to be equal to $q(r+1)$, attained at $b=-1$.
Now let $\beta>b^2$ and $X\succeq 0$ is equivalent to 
$Z:={\beta\over \beta-b^2} X_p - {b^2\over \beta-b^2}J\in\EE_p$.
As above, Lemma \ref{econdition} permits to bound the entries of $Z$  as follows:
${\beta\over \beta-b^2}c_s -{b^2\over \beta-b^2}\ge -\gamma$ for $1\le s\le q/2$.
Therefore,  the program (\ref{sdp1}) is upper bounded by
 \begin{equation}\label{sdpcw3}  
\begin{array}{ll}
 \underset{b,c,\beta}{\max}    & \frac{q}{2}(1-q\beta ) - pqb -c pq/2\\
\text{s.t. } & \beta(c+\gamma) \ge b^2(\gamma+1)\\
& b^2 < \beta\le 1, \ -1 \le b,c\le 1.
\end{array}
\end{equation}
At optimality, equality $\beta(c+\gamma)=b^2(\gamma+1)$ holds. This permits to express $c$ in terms of $b,\beta$ and to rewrite the objective function of (\ref{sdpcw3}) as ${q\over 2}(1-q\beta)-{pq\over 2}(b^2{\gamma+1\over \beta} +2b-\gamma)$.
For fixed $\beta$, the maximum of this quadratic function in $b$ is attained at $b=-{\beta\over \gamma +1}\in [-1,1]$ and is equal to 
${q\over 2}(1-q\beta)+{pq\over 2}({\beta\over \gamma +1} +\gamma)= {q\over 2} (\beta({p\over \gamma+1}-q)+p\gamma +1)$. As $q\le 2r$, we have ${p\over \gamma+1}\le q$ and thus  the latter quantity is maximized when $\beta=1$, so that the maximum of (\ref{sdpcw3}) is equal to 
${pq\over 2}(\gamma+{1\over \gamma+1}) -{q(q-1)\over 2}$. Hence, 
using $q\le 2r$ and $\gamma +{1\over \gamma+1}\le {3\over 2}$, we deduce that 
 this maximum is upper bounded by $3q(r+1)$. This concludes the proof  that the integrality gap of the clique-web inequality  is at most 3 when $q$ is even.

Consider now the case when q is odd. Then  $p$ is even and  $\text{W}^r_p$ consists of  $(q-1)/2$  bands of size $p$ and one  band of size $p/2$. The treatment is analogous to the case $q$ even, except we must replace the objective function in (\ref{obj1}) by
\begin{equation}
\label{obj2}
{q\over 2}(1-q\beta)-pqb-p(c_1+\ldots + c_{q-1\over 2})-{p\over 2}c_{q+1\over 2}
\end{equation}
and, as $p$ is even, the values on the bands can only be claimed to lie in $[-1,1]$ by Lemma \ref{econdition} (which amounts to setting $\gamma=1$ in the above argument). Specifically, 
if $\beta=0$  we can upper bound the objective function (\ref{obj2})  by  $2q(r+1)$ and, 
 if $\beta=b^2$, we can upper bound (\ref{obj2}) by $q(r+1)$.  Finally, if $\beta>b^2$, as above we do a Schur complement  and obtain
 ${\beta\over \beta-b^2}c_s-{b^2\over \beta-b^2} \ge -1$, so that 
 (\ref{obj2}) is upper bounded by the program (\ref{sdpcw3}) setting there $\gamma=1$.
 Hence the integrality gap of the clique-web inequality is also bounded by 3 for $q$ odd.
\qed\end{pf}

We conclude with several  remarks. \\
$\bullet$ We  just showed:
(i) $\ka(\text{\rm CW}^p_r,-e)\le 1+{(2r+1)^2\over q(2r+2)}$, and
(ii) $\ka(\text{\rm CW}^p_r,-e)\le {q+6r+5\over 4(r+1)}$ when $q\le 2r$.
Therefore,  asymptotically,  the integrality gap tends  to 1 as $q\rightarrow \infty$ and $r$ is fixed (by (i)), and it tends to $3/2$ as  $r\rightarrow \infty$ and $q$ is fixed (by (ii)).\\
$\bullet$ Our analysis is tight for  $q=2$, the case  of  bicycle odd wheel inequalities (since then the program (\ref{sdpcw3}) is equivalent to (\ref{sdp1}));  that is,
$ \sdp(\text{\rm CW}^r_{2r+3},-e)= -1 + p (\cos(\pi/p)+ {1\over \cos(\pi/p)+1} )$ (as mentioned  in \cite{DP93b}). This explains why in the proof of Theorem \ref{cwthm}    we use the precise estimate $-\cos\frac{\pi}{p} \le c,$ as opposed to the trivial bound $-1 \le c$, which was equally good for the our asymptotic analysis.  \\
$\bullet$ Pitowsky \cite{Pi08} shows the asymptotic lower bound $4/\pi\sim 1.27$ for the integrality gap of the clique-web inequality with $q=2r$ and $r\rightarrow \infty$.\\
$\bullet$ Given $b\in\oZ^n$ with $\sum_ib_i=1$, the  {\em hypermetric inequality} $w^Tx:= -\displaystyle\sum_{1\le i<j\le n}b_ib_jx_{ij}\le (\sum_i b_i^2-1)/ 2$ satisfies $\sdp(K_n,w) \le \sum_ib_i^2/2$;  thus its integrality gap is at most $3/2$, with equality if and only if $b=(1,1,-1,0,\ldots,0)$ (the case of triangle inequalities).



\begin{thebibliography}{1}

\bibitem{AMMN}
N.~Alon, K.~Makarychev, Y.~Makarychev, and A.~Naor.
\newblock Quadratic forms on graphs.
\newblock {\em Inventiones Mathimaticae}, 163(3):499--522, 2006.

\bibitem{AN04}
N.~Alon and A.~Naor.
\newblock Approximating the cut-norm via Grothendieck's inequality.
\newblock In {\em STOC}, 
pages 72--80,
  2004.

\bibitem{AB05}
S.~Arora, E.~Berger, G.~Kindler, M.~Safra, and E.~Hazan.
\newblock On non-approximability for quadratic programs.
\newblock In {\em FOCS}, pages 206--215, 2005.

\bibitem{B83}
F.~Barahona.
\newblock The max-cut problem on graphs not contractible to $K_5$.
\newblock {\em Oper. Res. Lett.}, 2(3):107--111, 1983.

\bibitem{BM}
F.~Barahona and A.~Mahjoub.
\newblock On the cut polytope.
\newblock {\em Math. Program.}, 36:157--173, 1986.

\bibitem{BJT93}
W.W. Barrett, C.R. Johnson, and P.~Tarazaga.
\newblock The real positive definite completion problem for a simple cycle.
\newblock {\em Linear Algebra Appl.}, 192:3--31, 1993.

\bibitem{BMMN}
M.~Braverman, K.~Makarychev, Y.~Makarychev, and A.~Naor.
\newblock The Grothendieck constant is strictly smaller than Krivine's bound.
Preprint, 	arXiv:1103.6161v2.

\bibitem{BOV}
J.~Bri\"et, F.~de Oliveira Filho, and F.~Vallentin.
Grothendieck inequalities for semidefinite programs with rank constraints.
Preprint, arXiv:1011.1754v1.

\bibitem{DL97}
M.M. Deza and M.~Laurent.
\newblock {\em Geometry of Cuts and Metrics}.
\newblock Springer, 1997.

\bibitem{FR94}
P.C. Fishburn and J.A. Reeds.
\newblock Bell inequalities, Grothendieck's constant, and root two.
{\em SIAM J. Disc. Math.,}7(1):48--56,1994.

\bibitem{GW95}
M.X. Goemans and D.P. Williamson.
\newblock Improved approximation algorithms for maximum cut and satisfiability
  problems using semidefinite programming.
  {\em  J. ACM,} 42:1115--1145, 1995.

\bibitem{Gro53}
A.~Grothendieck.
\newblock R\'esum\'e de la th\'eorie m\'etrique des produits tensoriels
  topologiques.
\newblock {\em Bol. Soc. Mat. Sao Paolo.}, 8:1--79, 1953.

\bibitem{Kri77}
J.~Krivine.
\newblock Sur la constante de Grothendieck.
\newblock {\em C. R. Acad. Sci. Paris S\'er. A-B,} 284(8):A445ÐA446,
1977.
  
\bibitem{Lau97}
M.~Laurent.
\newblock The real positive semidefinite completion problem for series-parallel
  graphs.
\newblock {\em Linear Algebra  Appl.}, 252:347--366, 1997.

\bibitem{Lau04}
M. Laurent.
Semidefinite relaxations for Max-Cut. In {\em The Sharpest Cut: The Impact of Manfred Padberg and His Work.} 
M. Gr\"otschel, ed., pages 257-290, MPS-SIAM Series in Optimization 4, 2004. 

\bibitem{lov79}
L.~Lov\'asz.
\newblock On the Shannon capacity of a graph.
\newblock {\em IEEE Trans. Inform. Th.}, IT-25:1--7, 1979.

\bibitem{Pi08}
I. Pitowsky.
New Bell inequalities for the singlet state: Going beyond the Grothendieck bound.
{\em Journal of Mathematical Physics}, 49:012101, 2008.

\bibitem{DP93c}
S.~Poljak. C.~Delorme.
\newblock Combinatorial properties and the complexity of a max-cut
  approximation.
\newblock {\em European J. Combin.}, 14:313--333, 1993.

\bibitem{DP93b}
S.~Poljak. C.~Delorme.
\newblock The performance of an eigenvalue bound on the max-cut problem in some
  classes of graphs.
\newblock {\em Discrete Math.}, 111:145--156, 1993.


\bibitem{Reeds}
J.A. Reeds.
\newblock A new lower bound on the real Grothendieck constant.
\newblock Preprint, 1991. Available at \verb1http://www.dtc.umn.edu/ reedsj/bound2.dvi1

\bibitem{W37}
K.~Wagner.
\newblock \"Uber eine Eigenschaft der ebenen Komplexe.
\newblock {\em Math. Ann.}, 114(1):570--590, 1937.

\bibitem{W06}
S.~Wehner.
\newblock Tsirelson bounds for generalized Clauser-Horne-Shimony-Holt inequalities.
\newblock {\em Phys. Rev. A}, 73:022110, 2006.

\end{thebibliography}


\end{document}